\documentclass{ifacconf}

\usepackage{graphicx}      
\usepackage{natbib}        
\usepackage{amsmath}
\usepackage{amssymb}
\usepackage{enumitem}
\usepackage{savesym}
\savesymbol{AND}
\usepackage{algorithmic}
\usepackage{algorithm}
\usepackage{tikz}
\usepackage{caption}
\usepackage{subcaption} 
\usepackage{setspace}

\let\Algorithm\algorithm
\renewcommand\algorithm[1][]{\Algorithm[#1]\setstretch{1.25}}

\usetikzlibrary{decorations.pathreplacing,shapes,arrows,calc}

\tikzstyle{block} = [draw, fill=white!20, rectangle, 
    minimum height=1cm, minimum width=1cm,rounded corners=5pt]
\tikzstyle{sum} = [draw, fill=white!20, circle, node distance=1cm]
\tikzstyle{input} = [coordinate]
\tikzstyle{output} = [coordinate]
\tikzstyle{tmp} = [coordinate]
\tikzstyle{pinstyle} = [pin edge={<-,line width=0.3mm,black}]

\newcommand{\algorithmiconline}{\textbf{Online: }}
\newcommand{\ONLINE}{\STATE \algorithmiconline}
\newcommand{\algorithmicoffline}{\textbf{Offline: }}
\newcommand{\OFFLINE}{\STATE \algorithmicoffline}

\begin{document}
\begin{frontmatter}

\title{Particle Model Predictive Control: Tractable Stochastic Nonlinear Output-Feedback MPC}  

\author{Martin A. Sehr \& Robert R. Bitmead} 
\address{Department of Mechanical \& Aerospace Engineering, University of California, San Diego, La Jolla, CA 92093-0411, USA \\
(e-mail: \{msehr, rbitmead\}@ ucsd.edu).}

\begin{abstract}                
We combine conditional state density construction with an extension of the Scenario Approach for stochastic Model Predictive Control to nonlinear systems to yield a novel particle-based formulation of stochastic nonlinear output-feedback Model Predictive Control. Conditional densities given noisy measurement data are propagated via the Particle Filter as an approximate implementation of the Bayesian Filter. This enables a particle-based representation of the conditional state density, or information state, which naturally merges with scenario generation from the current system state. This approach attempts to address the computational tractability questions of general nonlinear stochastic optimal control. The Particle Filter and the Scenario Approach are shown to be fully compatible and -- based on the time- and measurement-update stages of the Particle Filter -- incorporated into the optimization over future control sequences. A numerical example is presented and examined for the dependence of solution and computational burden on the sampling configurations of the densities, scenario generation and the optimization horizon.
\end{abstract}


\end{frontmatter}

\section{Introduction}
Model Predictive Control (MPC), in its original formulation, is a full-state feedback law (see \cite{mayne2000constrained,mayne2014model,precon}). This underpins two theoretical limitations of MPC: accommodation of output-feedback, and extension to include a compelling robustness theory given the state dimension is fixed. This paper addresses the first of these issues in a rather general, though practical setup. 

There has been a number of approaches to output-feedback MPC, mostly hinging on the replacement of the measured true state by a state estimate, which is computed via Kalman filtering (e.g. \cite{sehr2016sumptus,yan2005incorporating}), moving-horizon estimator (e.g. \cite{copp2014nonlinear,sui2008robust}), tube-based minimax estimators (e.g. \cite{mayne2009robust}), etc. Apart from \cite{copp2014nonlinear}, these designs, often for linear systems, separate the estimator design from the control design. The control problem may be altered to accommodate the state estimation error by methods such as: constraint tightening as in \cite{yan2005incorporating}, chance/probabilistic constraints  as in \cite{cannon2012stochastic} or \cite{schwarm1999chance}, and so forth. Likewise, for nonlinear problems, where the state estimation behavior is affected by control signal properties, the control may be modified to enhance the excitation properties of the estimator, as suggested in \cite{chisci2001systems,marafioti2014persistently}. Each of these aspects of accommodation is made in an isolated fashion.

The stochastic nonlinear output-feedback MPC algorithm presented in this paper is motivated by the structure of Stochastic Model Predictive Control (SMPC) via finite-horizon stochastic optimal control. The latter method requires propagating conditional state densities using a Bayesian Filter (BF) and solution of the Stochastic Dynamic Programming Equation (SDPE). By virtue of implementing a truly optimal finite-horizon control law in a receding horizon fashion, one can deduce a number of properties of the closed-loop dynamics, including recursive feasibility of the SMPC controller, stochastic stability and bounds characterizing closed-loop infinite-horizon performance, as discussed in \cite{sehr2016stochastic}. 

Unfortunately, solving for the stochastic optimal output-feedback controller, even on the finite horizon, is computationally intractable except for special cases such as linear quadratic Gaussian MPC because of the need to solve the SDPE, which incorporates the duality of the optimal control law in its effect on state observability. While the BF, required to propagate the conditional state densities, is readily approximated using a Particle Filter (PF), open-loop solution of the SDPE results in the loss of the duality of the optimal control. While not discussed in this paper, this effect can be mitigated sub-optimally by imposing excitation requirements as in \cite{chisci2001systems,marafioti2014persistently}.

Approximately propagating the conditional state densities by means of the PF naturally invites combination with the more recent advances in Scenario Model Predictive Control (SCMPC), as discussed for instance by \cite{blackmore2010probabilistic,calafiore2013stochastic,grammatico2016scenario,lee2014robust,mesbah2014stochastic,schildbach2014scenario}. Scenario methods deal with optimization of difficult, non-convex problems in which the initial task is recast as a parametrized collection of simpler, generally convex problems. Random sampling of uncertain signals and parameters is performed and the resulting collection of deterministic problem instances is solved. The focus has been on full state feedback for systems with linear dynamics and probabilistic state constraints. The technical construction is to take a sufficient number of samples (scenarios) to provide an adequate reconstruction of future controlled state densities for design. 

In contrast to solving the SDPE underlying the stochastic optimal control problem, the future controlled state densities in SCMPC are open-loop constructions. However, they present a natural fit combined with the particle-based conditional density approximations generated by the PF, where individual particles can be interpreted as scenarios from an estimation perspective. Moreover, while SCMPC is typically formulated in the linear case, the basic idea extends to the nonlinear case, albeit with the loss of many computation-saving features. In this paper, we propose and discuss this output-feedback version of SCMPC combined with the PF, which we call Particle Model Predictive Control (PMPC). Compared with the stochastic optimal output-feedback controller (computed via BF and SDPE), the PMPC controller is suboptimal in not accommodating future measurement updates and thereby losing both exact constraint violation probabilities along the horizon and the probing requirement inherent to stochastic optimal control. On the other hand, PMPC enables a generally applicable and, at least for small state dimensions, computationally tractable alternative for nonlinear stochastic output-feedback control.

The structure of the paper is as follows. We briefly introduce the problem setup in Section~\ref{sec:ps} and SMPC in Section~\ref{sec:SMPC} and proceed by introducing the PMPC control algorithm based on its individual components and parameters in Section~\ref{sec:PMPC}. After describing the algorithm and its correspondence to SMPC, we use a challenging scalar nonlinear example to demonstrate computational tractability and dependence of the proposed PMPC closed-loop behavior on a number of parameters in Section~\ref{sec:eg}. The example features nonlinear state and measurement equations and probabilistic state constraints under significant measurement noise. Finally, we conclude with Section~\ref{sec:conclusions}.

\section{Stochastic Optimal Control -- Setup}
\label{sec:ps}
We consider receding horizon output-feedback control for nonlinear stochastic systems of the form
\begin{align}
x_{t+1}&=f(x_t,u_t,w_t),\quad x_0\in\mathbb{R}^n,\label{eq:state}\\
y_t&=h(x_t,v_t), \label{eq:output}
\end{align}
starting from known initial state probability density function, $\pi_{0|-1} = \operatorname{pdf}(x_0)$. To this end, we denote the data available at time $t$ by 
\begin{align*}
\mathbf{\zeta}^t&\triangleq\{y_0,u_0,y_1,u_1,\dots,u_{t-1},y_t\},&\mathbf{\zeta}^0&\triangleq\{y_0\}.
\end{align*}
The \textit{information state,} denoted $\pi_t$, is the conditional density of state $x_t$ given data $\mathbf{\zeta}^t$.
\begin{align}\label{eq:pikk}
\pi_{t}&\triangleq\operatorname{pdf}\left(x_{t}\mid \mathbf{\zeta}^t \right).
\end{align}
We further impose the following standing assumption on the random variables and control inputs.
\begin{assum}\label{assm:sys}
The signals in~(\ref{eq:state}-\ref{eq:output}) satisfy:
\begin{enumerate}[label=\arabic*.]
\item $\{w_t\}$ and $\{v_t\}$ are sequences of independent and identically distributed random variables.
\item $x_0, w_t, v_l$ are mutually independent for all $t,l\geq 0$.
\item The control input $u_t$ at time instant $t\geq 0$ is a function of the data $\mathbf{\zeta}^t$ and given initial state density $\pi_{0\mid -1}$.
\end{enumerate}
\end{assum}
Denote by $\mathbb{E}_t[\,\cdot\,]$ and $\mathbb{P}_t[\,\cdot\,]$ the conditional expected value and probability with respect to state $x_t$ -- with conditional density $\pi_t$ -- and random variables $\{(w_k,v_{k+1}):k\geq t\}$, respectively, and by $\epsilon_{k}$ the constraint violation level of constraint $x_k \in \mathbb{X}_k$. Our goal is to solve the \emph{finite-horizon stochastic optimal control problem} (FHSOCP)
\begin{multline*}
\!\!\!\!\mathcal{P}_{N}(\pi_{t}): 
\left\{\!\!\begin{array}{cl}
\inf_{u_{t},\ldots,u_{t+N-1}}&\!\! \mathbb{E}_t\left[\sum_{k=t}^{t+N-1}{c(x_{k},u_{k})} + c_{N}(x_{t+N})\right],\\[0.2cm]
\text{s.t.}&\!\! x_{k+1}=f(x_{k},u_{k},w_{k}),\\[0.1cm]
&\!\! x_{t}\sim\pi_t,\\[0.1cm]
 &\!\! \mathbb{P}_{k+1}\left[ x_{k+1} \in \mathbb{X}_{k+1} \right]\geq 1 - \epsilon_{k+1},\\[0.1cm]
   &\!\! u_{k}\in\mathbb{U}_{k},\\[0.1cm]
&\!\! k=t,\dots,t+N-1.
\end{array}\right.
\end{multline*}
In theory, solving the FHSOCP at each time $t$ and subsequently implementing the first control in a receding horizon fashion leads to a number of desirable closed-loop properties, as discussed in \cite{sehr2016stochastic}. However, solving the FHSOCP is computationally intractable in practice, a fact that has led to a number of approaches in MPC for nonlinear stochastic dynamics. We propose a novel strategy that is oriented at the structure of SMPC based on the FHSOCP, but numerically tractable at least for low state dimensions.

As a result of the Markovian state equation \eqref{eq:state} and measurement equation \eqref{eq:output}, the optimal control inputs in the FHSOCP must inherently be \textit{separated} feedback policies (e.g.~\cite{bertsekas1995dynamic,BKKUM1986}). That is, control input $u_{t}$ depends on the available data $\mathbf{\zeta}^t$ and initial density $\pi_{0\mid -1}$ solely through the current information state, $\pi_{t}$. Optimality thus requires propagating $\pi_{t}$ and policies $g_t$, where
\begin{align}\label{eq:policies}
u_t = g_t(\pi_{t}).
\end{align}
Motivated by this two-component separated structure of stochastic optimal output-feedback control, we propose an extension of the SCMPC approach to nonlinear systems, merged with a numerical approximation of the information state update via particle filtering. Before proceeding with this novel approach, we briefly revisit the two components of SMPC via solution of the FHSOCP.

\section{Stochastic Model Predictive Control}
\label{sec:SMPC}
The information state is propagated via the \emph{Bayesian Filter} (see e.g.~\cite{chen2003bayesian,simon2006optimal}):
\begin{align}
\pi_{t}&=
\frac{\operatorname{pdf}(y_{t}\mid x_{t})\,\pi_{t\mid t-1}}{\int \operatorname{pdf}(y_{t}\mid x_{t})\,\pi_{t\mid t-1}\,dx_{t}},\label{eq:BF_rec} \\
\pi_{t+1\mid t} &\triangleq
\int \operatorname{pdf}(x_{t+1} \mid x_{t},u_{t}) \,\pi_{t}\, dx_{t},\label{eq:BF_pred}
\end{align}
for $t\in\{0,1,2,\ldots\}$ and initial density $\pi_{0\mid -1}$.
The recursion~(\ref{eq:BF_rec}-\ref{eq:BF_pred}) has the following features:
\begin{itemize}
\item The \emph{measurement update}~\eqref{eq:BF_rec} combines the \emph{a priori} conditional density, $\pi_{t\mid t-1}$, and $\operatorname{pdf}(y_{t}\mid x_{t})$, derived from~\eqref{eq:output} using knowledge of: the function $h(\cdot,\cdot)$, the density of $v_{t}$, and the value of $y_t$. 
\item The \emph{time update}~\eqref{eq:BF_pred} combines $\pi_{t}$ and $\operatorname{pdf}(x_{t+1} \vert x_{t},u_{t})$, derived from~\eqref{eq:state} using knowledge of: control input $u_{t}$, function $f(\cdot,\cdot,\cdot)$, and the density of $w_t$. 
\item For linear Gaussian systems, the filter recursion~(\ref{eq:BF_rec}-\ref{eq:BF_pred}) reduces to the well-known Kalman Filter.
\end{itemize}
Combined with solution of the FHSOCP, this leads to the following SMPC algorithm, as 
discussed in \cite{sehr2016stochastic}.

\begin{algorithm}[H] \label{algo:SMPC}                   
\caption{Stochastic Model Predictive Control}\label{RHSOC}                       
\begin{algorithmic}[1]
\OFFLINE
\STATE 
Solve $\mathcal{P}_N(\cdot)$ for the first optimal policy, $g_0^{\star}(\cdot)$.
\ONLINE
\FOR{$t=0,1,2,\ldots$}
\STATE Measure $y_t$
\STATE Compute $\pi_t$
\STATE Apply first optimal control policy, $u_{t} = g_0^{\star}(\pi_{t})$
\STATE Compute $\pi_{t+1\mid t}$
\ENDFOR
\end{algorithmic}
\end{algorithm}
Notice how this algorithm differs from common practice in stochastic model predictive control in that it explicitly uses the information states $\pi_t$. Throughout the literature, these information states -- conditional densities -- are commonly replaced by state estimates. While this makes the problem more tractable, one no longer solves the underlying stochastic optimal control problem.
The central divergence however lies in Step~2 of the algorithm, in which the SDPE is presumed solved offline for the optimal feedback policies, $g_t(\pi_t)$, from \eqref{eq:policies}. This is an extraordinarily difficult proposition in many cases but captures the optimality, and hence duality, as a closed-loop feedback control law. The complexity of this step lies not only in computing a vector functional but also in the internal propagation of the information state within the SDPE.

\section{Tractable Nonlinear Output-Feedback Model Predictive Control}
\label{sec:PMPC}
In this section, we motivate a novel approach to output-feedback MPC that maintains the separated structure of SMPC while being numerically tractable for modest problem size.

\subsection{Approximate Information State \& Particle Filter}
The BF~(\ref{eq:BF_rec}-\ref{eq:BF_pred}) propagates the information state $\pi_t$ to implement a necessarily separated stochastic optimal output-feedback control law. While implementing this recursion precisely is possible only in special cases such as linear Gaussian systems, where the densities can be finitely parametrized, the BF can be implemented approximately by means of the Particle Filter, with the approximation improving with the number of particles, as described for instance in \cite{simon2006optimal}. In parallel with the BF, the PF consists of two parts: the forward propagation of the state density, and the resampling of the density using the next measurement. 

The following algorithm describes a version of the PF amenable to PMPC in the context of this paper. This is a slightly modified version of the filter design described by \cite{simon2006optimal}.
\begin{algorithm}[H]                   
\caption{Particle Filter (PF)}                       
\begin{algorithmic}[1]
\STATE Sample $N_p$ particles, $\{x_{0,p}^-,\,p=1,\dots,N_p\}$, from density $\pi_{0\mid -1}$.
\FOR{$t=0,1,2,\ldots$}
\STATE Measure $y_k$.
\STATE Compute the relative likelihood $q_p$ of each particle $x_{t,p}^-$ conditioned on the measurement $y_t$ by evaluating $\operatorname{pdf}(y_t \mid x_{t,p}^-)$ based on~\eqref{eq:output} and $\operatorname{pdf}(v_t)$.
\STATE Normalize $q_p \to q_p/\sum_{p=1}^{N_p}{q_p}$.
\STATE Sample $N_p$ particles, $x_{t,p}^+$, via \emph{resampling} based on the relative likelihoods $q_p$.
\STATE Given $u_t$, propagate $x_{t+1,p}^- = f(x_{t,p}^+, u_t, w_{t,p})$, where $w_{t,p}$ is generated based on $\operatorname{pdf}(w_t)$.
\ENDFOR
\end{algorithmic}
\end{algorithm}
While a number of variations -- such as roughening of the particles and differing \emph{resampling} strategies, including importance sampling -- of  this basic algorithm may be sensible depending on the system at hand, this basic algorithm suffices in presenting a numerical method of approximating the Bayesian Filter to arbitrary degree of accuracy with increasing number of particles $N_p$ (see e.g. \cite{smith1992bayesian}). For a more detailed discussion on the PF for use in state-estimate feedback control, see \cite{rawlings2009model}.

\subsection{Scenario MPC and Particle Model Predictive Control}
\begin{figure}[htb]
  \centering
  \includegraphics[width=\columnwidth]{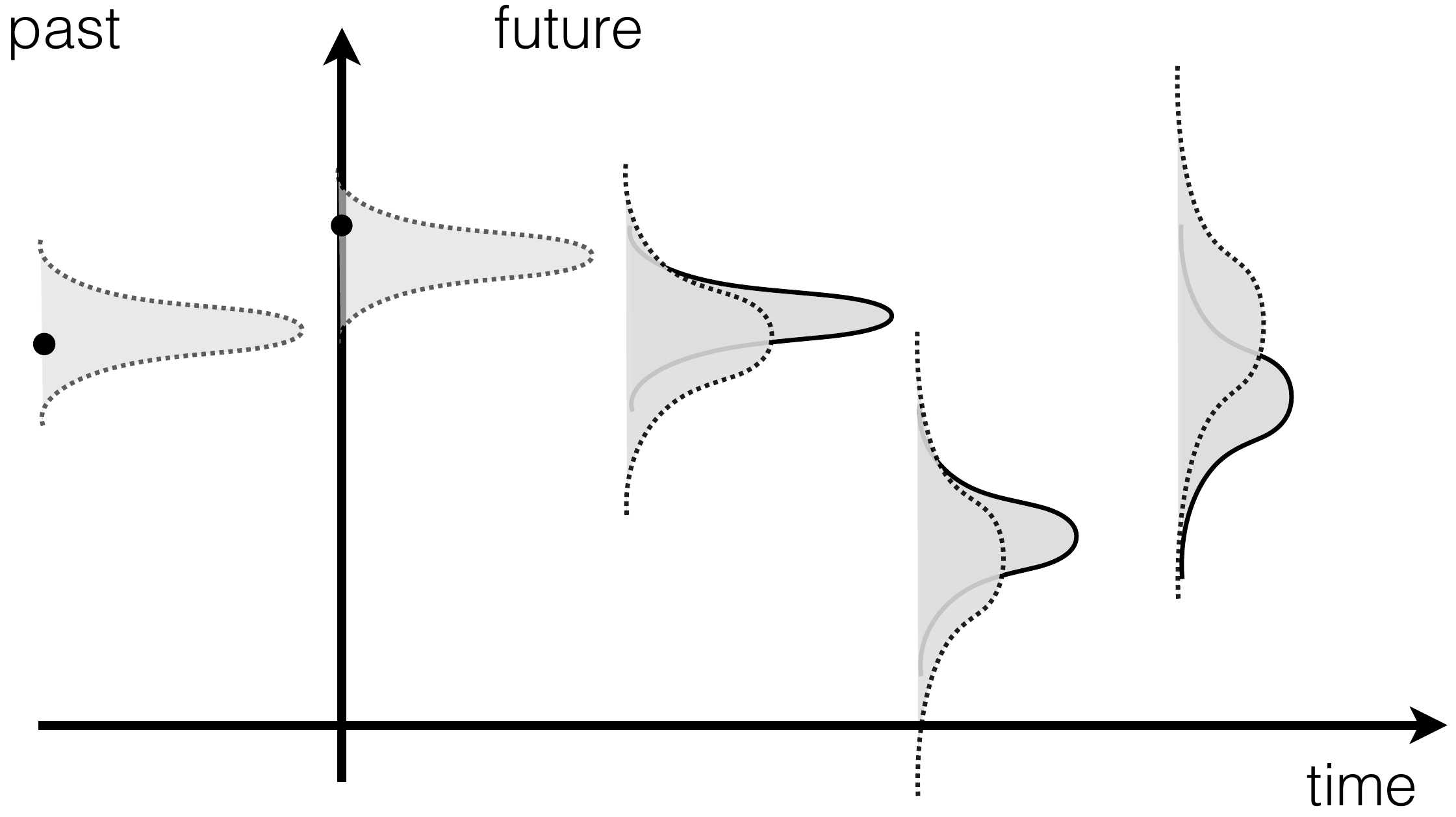}
  \caption{State density evolution in: Scenario MPC calculations (dots and solid outlines) and, Particle MPC (dashed outlines), for three steps into the future. \label{fig:PMPC}}
\end{figure}

The Scenario Approach to MPC (e.g. \cite{calafiore2013stochastic,grammatico2016scenario,lee2014robust,mesbah2014stochastic,schildbach2014scenario}) commences from state $x_t$ or state estimate, $\hat x_{t|t}$. It propagates, i.e. simulates, an open-loop controlled stochastic system with sampled process noise density $\operatorname{pdf}(w_t)$. These propagated samples are then used to evaluate controls for constraint satisfaction and for open-loop optimality with probabilities tied to the sampled $w_t$ densities. In many regards, this is congruent to repeated forward propagation of the PF via \eqref{eq:BF_pred} without measurement update \eqref{eq:BF_rec} and commencing from a singular density at $x_t$ or $\hat x_{t|t}$. Particle MPC simply replaces the starting point, $\hat x_{t|t},$ by the collection of particles $\{x^+_{t,p}, p=1,\dots.N_p\}$ distributed as $\pi_t$, as illustrated in Figure~\ref{fig:PMPC}.

Before introducing the PMPC algorithm, we define a sampled, particle version of the FHSOCP, with $N_s$ scenarios and $N_p$ available \emph{a posteriori} particles at time $t$,
\begin{multline*}
\!\!\!\!\tilde{\mathcal{P}}_{N}(\{x_{t,p}^+, p = 1,\ldots,N_p\}): \\
\left\{\!\!\begin{array}{cl}
\inf_{u_{t},\ldots,u_{t+N-1}}&\!\!
\sum_{s=0}^{N_s} \left(\sum_{k=t}^{t+N-1}{c(x_{k,s},u_{k})} + c_{N}(x_{t+N,s})\right),\\[0.2cm]
\text{s.t.}&\!\! x_{k+1,s}=f(x_{k,s},u_{k},w_{k,s}),\\[0.1cm]
&\!\! x_{t,s}\in\{x_{t,p}^+, p = 1,\ldots,N_p\},\\[0.1cm]
 &\!\! \tilde{\mathbb{P}}_{k+1}\left[ x_{k+1} \in \mathbb{X}_{k+1} \right] 
\geq 1 - \epsilon_{k+1},\\[0.1cm]
   &\!\! u_{k}\in\mathbb{U}_{k},\\[0.1cm] 
&\!\! s=1,\ldots,N_s,\quad  k=t,\dots,t+N-1,
\end{array}\right.
\end{multline*}
where the statement 
\begin{align*}
\tilde{\mathbb{P}}_{k+1}\left[ x_{k+1} \in \mathbb{X}_{k+1} \right]\geq 1 - \epsilon_{k+1}
\end{align*} 
means that $x_{k+1,s} \in \mathbb{X}_{k+1}$ for at least $(1 - \epsilon_{k+1})N_s$ scenarios. Following the approach in \cite{schildbach2013randomized}, one may also choose to replace this constraint by $x_{k+1} \in \mathbb{X}_{k+1}$ and select the number of scenarios $N_s$ according to the desired constraint violation levels $\epsilon_{k+1}$. We are now in position to formulate the PMPC algorithm following the schematic in Figure~\ref{fig:PMPC}.
\begin{algorithm}[H]                   
\caption{Particle Model Predictive Control (PMPC)}\label{PMPC}                       
\begin{algorithmic}[1]
\STATE Generate $N_p$ \emph{a priori} particles, $x_{0,p}^-$, based on $\pi_{0\mid -1}$.
\FOR{$t=0,1,2,\ldots$}
\STATE Measure $y_t$.
\STATE Compute the relative likelihood $q_p$ of each particle $x_{t,p}^-$ conditioned on the measurement $y_t$ by evaluating $\operatorname{pdf}(y_t \mid x_{t,p}^-)$ based on~\eqref{eq:output} and $\operatorname{pdf}(v_t)$.
\STATE Normalize $q_p \to q_p/\sum_{p=1}^{N_p}{q_p}$.
\STATE Generate $N_p$ \emph{a posteriori} particles, $x_{t,p}^+$, via \emph{resampling} based on the relative likelihoods $q_p$.
\STATE Solve $\tilde{\mathcal{P}}_{N}(\{x_{t,p}^+, p = 1,\ldots,N_p\})$ for the optimal scenario control values $u_t^{\star},\ldots,u_{t+N-1}^{\star}$.
\STATE Given $u_t^{\star}$, propagate $x_{t+1,p}^- = f(x_{t,p}^+, u_t^{\star}, w_{t,p})$, where $w_{t,p}$ is generated based on $\operatorname{pdf}(w_t)$.
\ENDFOR
\end{algorithmic}
\end{algorithm}

\subsection{Computational Demand}
Computational tractability of PMPC deteriorates with increasing: number of particles; number of scenarios; system dimensions; control signal grid spacing; MPC horizon. While the number of particles required for satisfactory performance of the PF grows exponentially with the state dimension (e.g.~\cite{snyder2008obstacles}), it is unclear how to select an appropriate number of scenarios in the nonlinear case. Suppose the state and input dimensions are $n$ and $m$ and that the numbers of particles and scenarios are chosen as $N_p = P^n$ and $N_s = S^n$ for positive integers $P$ and $S$, respectively, and that the MPC horizon is $N$. Further assuming a grid of $U^m$ points in the control space and brute-force evaluation of all possible sequences, the order of growth for PMPC is approximately $\mathcal{O}(P^n + S^nU^{mN})$. 

Notice that the computational demand associated with the conditional density approximation in PMPC is additive in terms of the overall computational demand. This indicates that, provided the PF is computationally tractable for given state dimensions, tractability of PMPC is roughly equivalent to tractability of standard state-feedback SCMPC. In the example below, we found that scenario optimization tends to be the computational bottleneck at least for low system dimensions. Clearly, this observation holds only when the scenario optimization is performed by explicit enumeration of all feasible sequences over a grid in the control space, which may be avoided for particular problem instances. But the experience also confirms that in the nonlinear case the open- or closed-loop control calculation dominates the computational burden in comparison to state estimation.

\section{Numerical Example}
\label{sec:eg}
\begin{figure*}[t!]
    \centering
    \begin{subfigure}[t]{0.5\textwidth}
        \centering
        \includegraphics[width=\columnwidth]{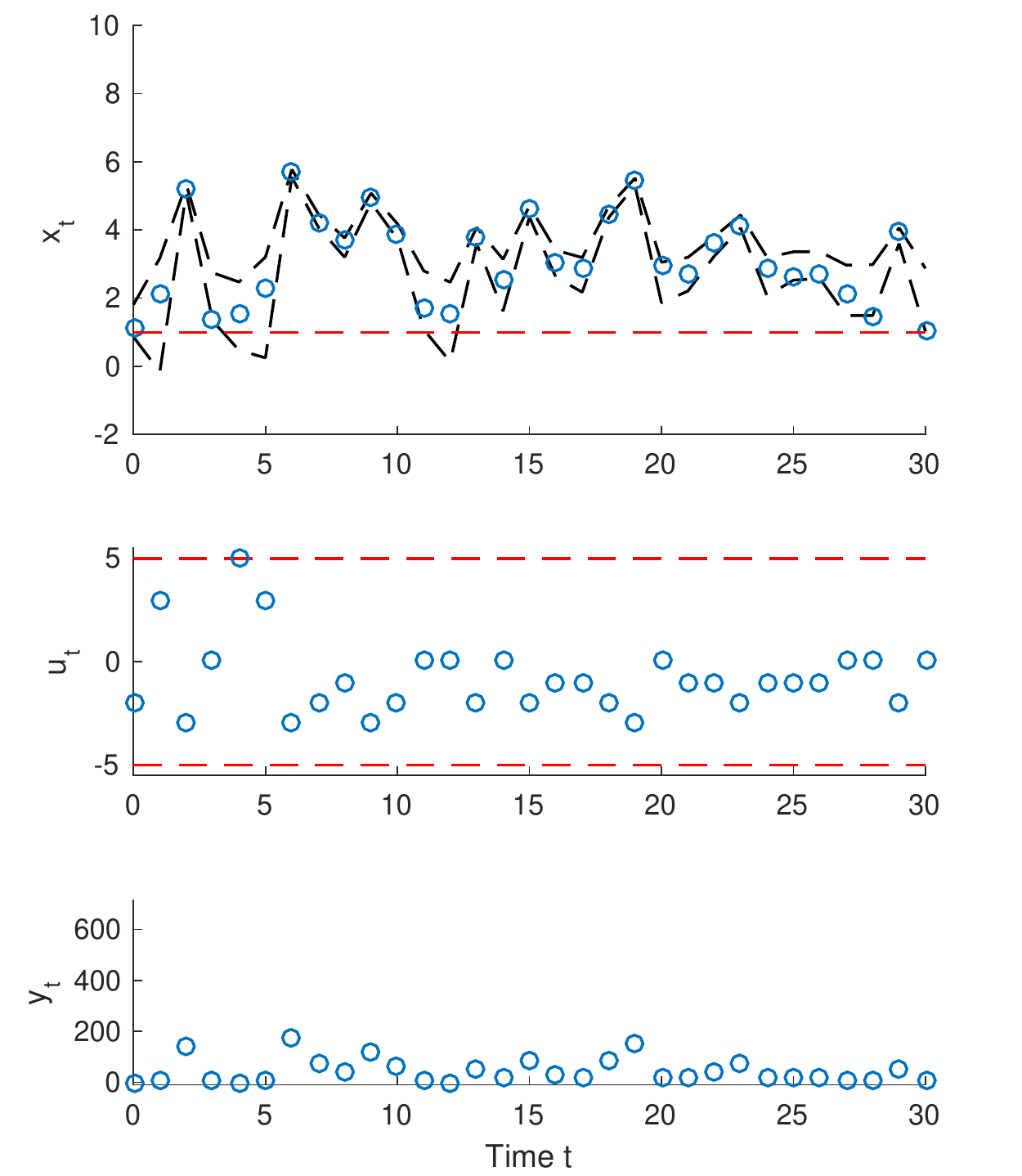}
        \caption{$N=3$, $N_p = 5,000$, $N_s = 1,000$.}
        \label{fig:densities1}
    \end{subfigure}%
    ~ 
    \begin{subfigure}[t]{0.5\textwidth}
        \centering
        \includegraphics[width=\columnwidth]{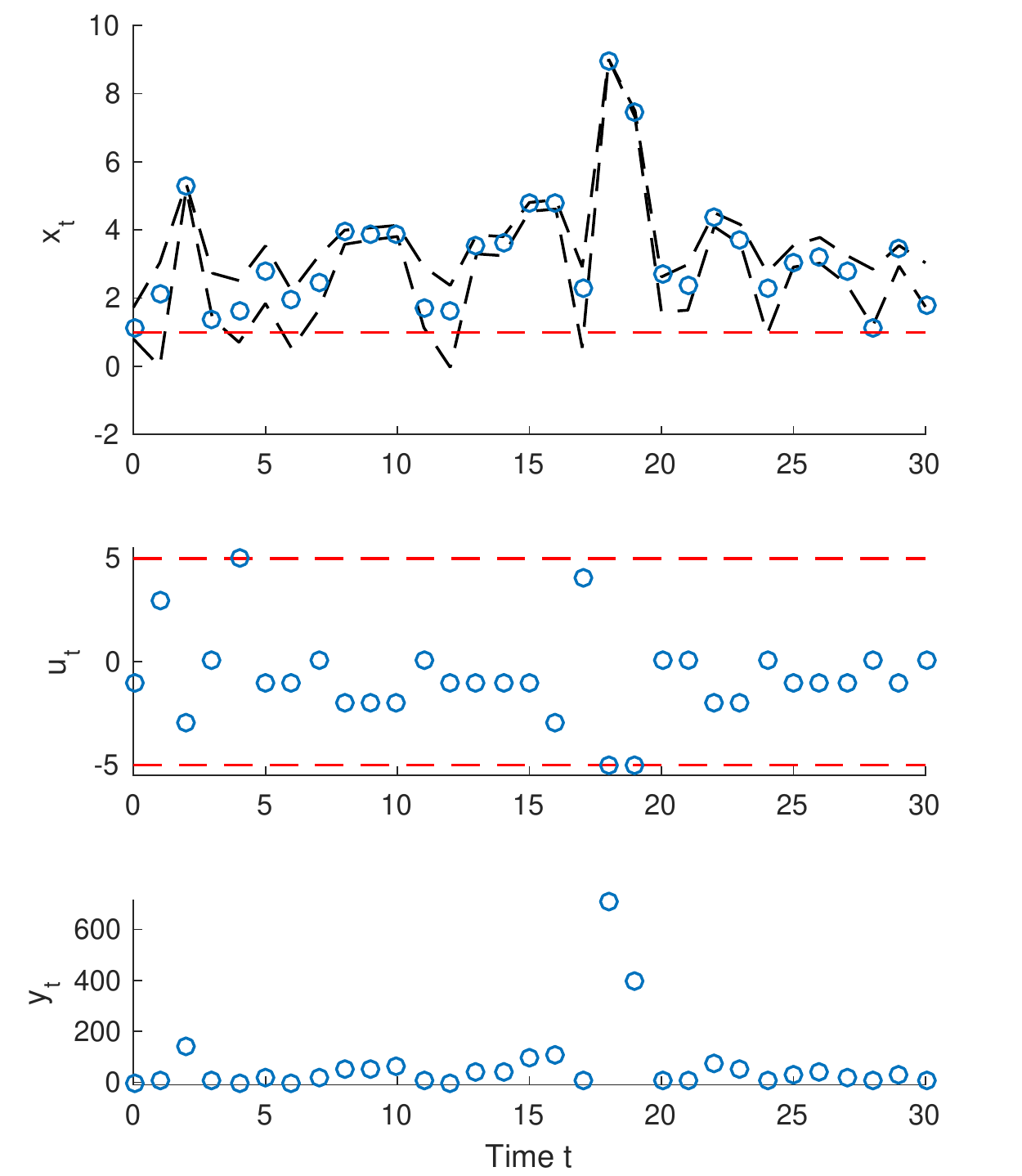}
        \caption{$N=3$, $N_p = 100$, $N_s = 1,000$.}
        \label{fig:densities2}
    \end{subfigure}
    \\
    \begin{subfigure}[t]{0.5\textwidth}
        \centering
        \includegraphics[width=\columnwidth]{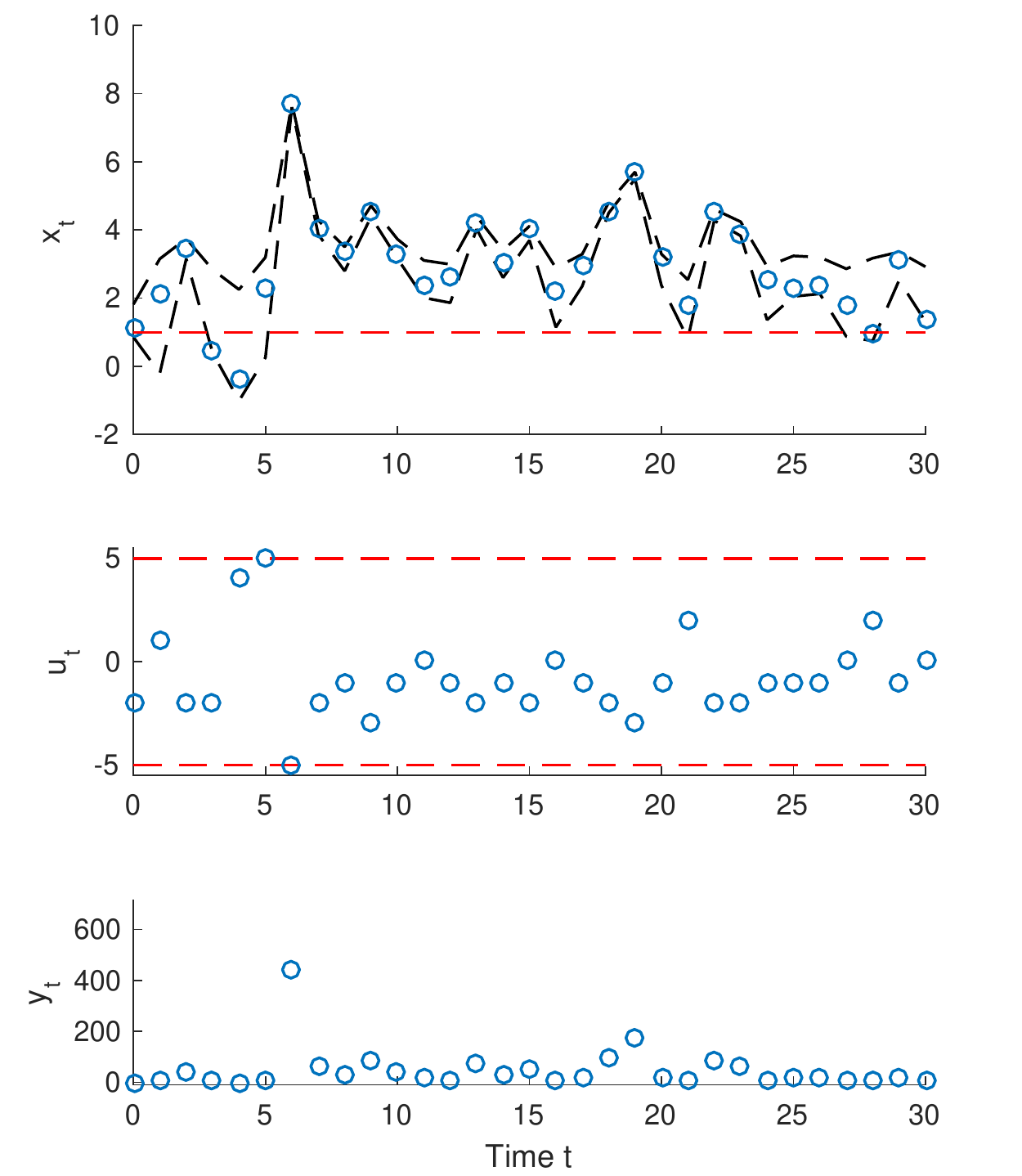}
        \caption{$N=3$, $N_p = 5,000$, $N_s = 50$.}
        \label{fig:densities3}
    \end{subfigure}%
    ~ 
    \begin{subfigure}[t]{0.5\textwidth}
        \centering
        \includegraphics[width=\columnwidth]{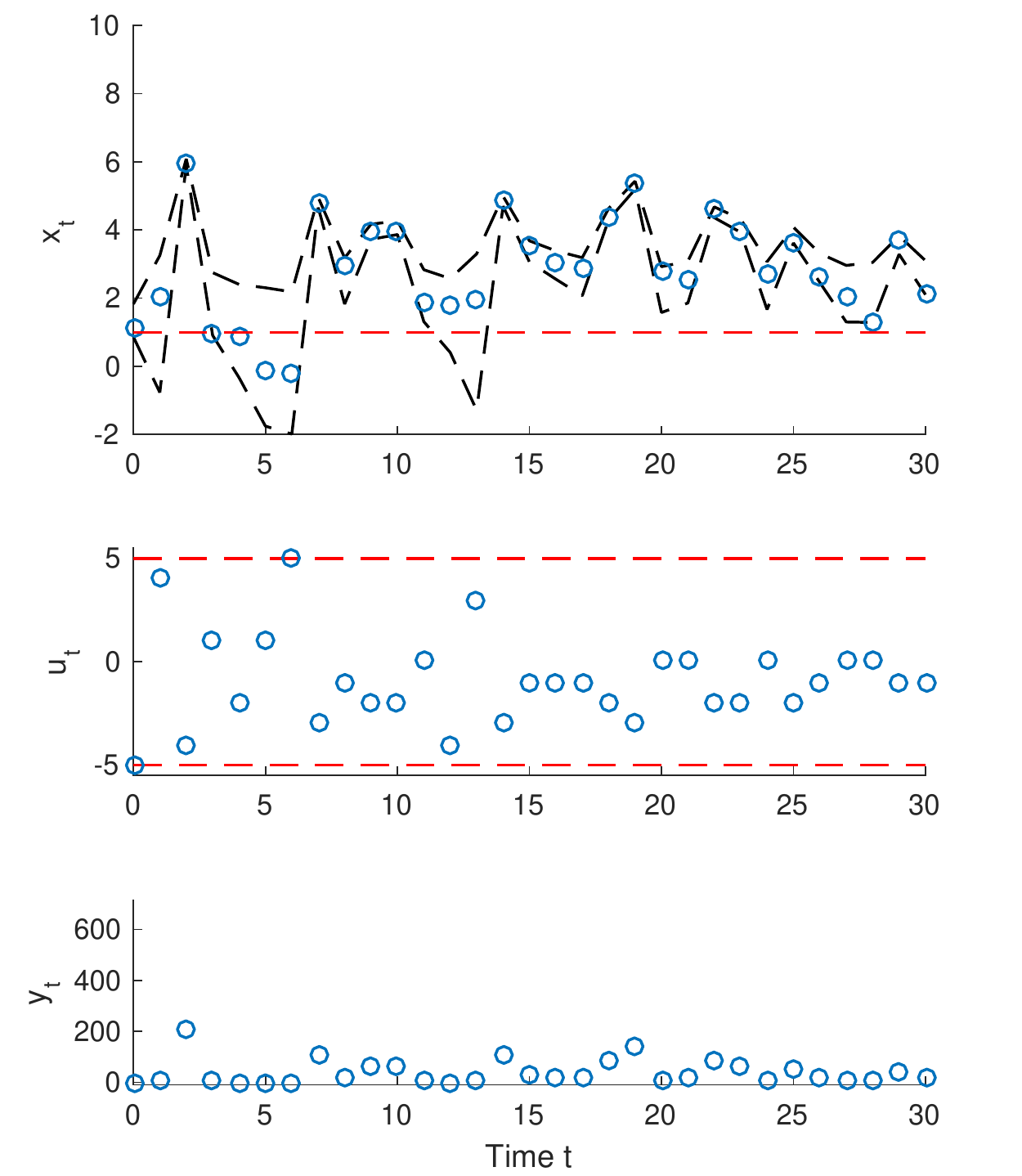}
        \caption{$N=2$, $N_p = 5,000$, $N_s = 1,000$.}
        \label{fig:densities4}
    \end{subfigure}
    \caption{Simulation data for example in Section~\ref{sec:eg} over $30$ samples, running PMPC with control horizon $N$, number of particles $N_p$ and number of scenarios $N_s$. State, control and measurement values (blue), probabilistic and hard constraints (red), $95\%$ confidence interval of PF (black). All controllers are subject to the same realization of the process noise $w_k$.}
\label{fig:simu}
\end{figure*}

Consider the scalar, nominally unstable nonlinear system
\begin{align*}
x_{t+1} &= 1.5\, x_t + \operatorname{atan}{\left((x_t - 1)^2\right)}\,u_t + w_t, \\
y_t &= x_t^3 - x_t + v_t,
\end{align*}
where $x_0,w_t$ and $v_l$ are mutually independent random variables for all $t,l \geq 0$ and
\begin{align*}
x_0 &\sim \mathcal{U}(1,2),&
w_t&\sim\mathcal{U}(-2,2),&
v_t &\sim \mathcal{N}(0,5),
\end{align*}
for all $t\geq 0$. We aim to minimize the quadratic cost function
\begin{multline*}
J_N(\pi_{t},u_t,\ldots,u_{t+N-1})= \\
\mathbb{E}_t\left[\sum_{k=t}^{t+N-1}{\left(100\, x_{k}^2 + u_{k}^2\right)} + 100\,x_{t+N}^2 \right],
\end{multline*}
while satisfying the constraints
\begin{align*}
\mathbb{P}_{k+1} [x_{k+1} \geq 1] &\geq 0.9, &
-5\leq u_{k} & \leq 5,
\end{align*}
along the control horizon $N$, that is $k\in\{t,\ldots,t+N-1\}$ for $t\geq 0$. Notice how this system has both limited observability and controllability close to the constraint but infeasible unconstrained optimal states. In combination with the very noisy measurements, this is a challenging control problem. To implement PMPC as described in Section~\ref{sec:PMPC} for this nonlinear stochastic output-feedback control problem, we further restrict the control inputs to integer values, such that $u_{t}\in\{-5,-4,\ldots,4,5\}$. Figure~\ref{fig:simu} displays simulated closed-loop state trajectories, control values and measurement values for four PMPC controllers with differing parameters subject to the same realizations of process and measurement noise, respectively.

Figure~\ref{fig:densities1} displays closed-loop simulation results under PMPC with horizon $N = 3$, $N_p = 5,000$ particles and $N_s = 1,000$ scenarios. While the poor observability properties of the system show close to the probabilistic constraint, it is satisfied at all times in this simulation. This is still the case when decreasing the number of particles to $N_p = 100$ in the simulation displayed Figure~\ref{fig:densities2}. However, we see how in this case, the decreased accuracy of the PF leads to larger state-values in closed-loop. Similar behavior is observed in Figure~\ref{fig:densities3} when reducing the number of scenarios to $N_s = 50$. Additionally, the controller violates the probabilistic constraint $3$ times in this case. This trend continues when reducing the horizon to $N=2$, as displayed in Figure~\ref{fig:densities4}.

\section{Conclusion}
\label{sec:conclusions}
We presented PMPC as a novel approach to output-feedback control of stochastic nonlinear systems. Generating scenarios not only from the distribution of the process noise but also from the particles of the Particle Filter, PMPC combines the benefits of the Particle Filter and Scenario MPC in a natural fit, allowing for a numerically tractable version of stochastic MPC with general nonlinear dynamics, cost and probabilistic constraints. Given a particular system instance, the algorithm and its properties may be adapted to exploit specific problem structure. Such extensions include: sub-optimal probing via additional constraints; scenario removal; provable closed-loop properties such as constraint satisfaction with specified confidence levels; optimization over parametrized policies.

\bibliography{IFAC2017}            


\end{document}